\documentclass[12pt]{article}
\usepackage{amsmath}
\usepackage{amssymb}
\newsymbol \blackbox 1004
\newcommand{\eh}{\hfill}\newlength{\sperr}

\newenvironment{proof}{{\settowidth{\sperr}{\bf\rm
Proof}%
\par\addvspace{0.3cm}\noindent\parbox[t]{1.3\sperr}
{\bf\rm P\eh r\eh o\eh o\eh f\eh }%
}}{\nopagebreak\mbox{}
$\blackbox$\par\addvspace{0.3cm}}

\def\b{\beta}
\def\g{\gamma}
\def\G{\Gamma}

\def\s{\sigma}
\def\la{\lambda}

\def\T{\Theta}

\def\vp{\varphi}
\def\ve{\varepsilon}
\def\wh{\widehat}
\def\wt{\widetilde}
\def\ov{\overline}
\def\p{\partial}

\def\BC{{\mathbb C}}
\def\BR{{\mathbb R}}

\def\ker{{\rm Ker\ }}
\newtheorem{Pa}{Paper}[section]
\newtheorem{Tm}[Pa]{{\bf Theorem}}

\newtheorem{Cy}[Pa]{{\bf Corollary}}
\newtheorem{Rk}[Pa]{{\bf Remark}}

\newtheorem{Pn}[Pa]{{\bf Proposition}}
\newcommand{\CC}
{{\mathchoice {\setbox0=\hbox{$\displaystyle\rm
C$}\hbox{\hbox
to0pt{\kern0.4\wd0\vrule height0.9\ht0\hss}\box0}}
{\setbox0=\hbox{$\textstyle\rm C$}\hbox{\hbox
to0pt{\kern0.4\wd0\vrule height0.9\ht0\hss}\box0}}
{\setbox0=\hbox{$\scriptstyle\rm C$}\hbox{\hbox
to0pt{\kern0.4\wd0\vrule height0.9\ht0\hss}\box0}}
{\setbox0=\hbox{$\scriptscriptstyle\rm C$}\hbox{\hbox
to0pt{\kern0.4\wd0\vrule height0.9\ht0\hss}\box0}}}}

\title{On  explicit inversion of a subclass of operators with   $D$-difference kernels
and Weyl theory of the corresponding canonical systems}

\author{A.L. Sakhnovich, A.A. Karelin, J. Seck-Tuoh-Mora, 
\\
G. Perez-Lechuga, 
M. Gonzalez-Hernandez}

\date{}
\begin{document}
\maketitle

\begin{abstract} 
Explicit inversion formulas for a subclass of integral operators with $D$-difference
kernels on a finite interval are obtained. A case of the positive operators is treated in greater
detail. An application to the inverse problem to recover canonical
system from a Weyl function is given.
\end{abstract}

{MSC(2000) Primary 34A55,  45Q05; Secondary 47B65, 47G10}

{\it Keywords: integral operator with difference kernel, operator with $D$-difference kernel, explicit inversion,
canonical system, inverse problem, Weyl function.}

\section{Introduction} \label{intro}
\setcounter{equation}{0}
Integral operators with difference kernels are important
in mathematics and applications and are actively used in  the study of  numerous homogenious
processes. The papers \cite{GoKr, Kr} on the inversion of the operators with difference kernels
on the semi-axis became classical. Various results and references on the operators with difference kernels
on a finite  interval or a system of intervals are given in \cite{SaLdiff, SaL1}. Interesting explicit
results on the inversion of the operators with exponential type difference kernels on a
finite interval one can find in \cite{BGK, GoKaSc}.

Operators with  $D$-difference kernels in $L^2_p(0,l)$,  which  we  shall treat, are bounded operators
of the form
\begin{eqnarray}\label{0.1}
&&S_l f=Sf=\frac{d}{dx} \int_0^l
s(x,t)f(t)dt,  \quad  s(x,t)=\big\{ s_{i j}(x,t) \big\}_{i,j=1}^p, 
\\  \label{0.2}    && s_{ij}(x,t)=s_{ij}(d_ix-d_j t), \quad s_{ij}(x) \in L^2(-d_j l, \, d_i l), 
\end{eqnarray}
where  $D=D^*={\mathrm{diag}}\{d_1, \, d_2, \, \ldots , \, d_p \}>0$
is a fixed $p \times p$ diagonal matrix. 
The notion of an operator with  a $D$-difference kernel is a natural generalization of
the operator with a difference kernel, i.e., of  the case  $D=I_p$,  where  $I_p$ is the $p \times p$
identity matrix.  The class of operators with $D$-difference kernels on a finite interval
includes the operators with difference kernels on systems of intervals,
which are important, for instance, in elasticity theory,
diffraction theory, and the theory of stable processes (see \cite{Ka} and  Chapter 6 in \cite{SaL1}).

Explicit inversion formulas for  an interesting subclass of  operators with  $D$-difference kernels
are obtained in Section \ref{Inver} of  this paper using the classical results on semiseparable operators. 
Note also that the inversion of semiseparable matrices and operators is another
interesting and actively developed theory, see \cite{GGK1, GK84} and bibliography in \cite{VVGM}.
Some further possible applications are connected with the paper \cite{KarLT}.

Operator identities for the operators  with  $D$-difference kernels are discussed in Section \ref{OpId}.

The case of positive and boundedly invertible operators with $D$-difference kernels
is treated further in Theorem \ref{TmPos} of Section \ref{InvPr}.
As an application, we solve explicitly in terms of Weyl functions an inverse problem for a subclass of canonical
systems. Some results from \cite{SaL2, SaL3} are developed further in this section too.

 We use the standard notations $\BC$ and $\BC_+$  for the complex plane and upper semi-plane, respectively.
  By $\{{\cal H}_1, \,{\cal H}_2\}$ we denote
   the class of the bounded linear operators acting from ${\cal H}_1$ into ${\cal H}_2$,
   and by $\s(\b)$ we denote the spectrum of $\b$.

\section{Inversion of  operators with $D$-difference kernels} \label{Inver}
\setcounter{equation}{0}  
Consider a self-adjoint operator with $D$-difference kernel
\begin{equation}\label{i1}
S=I+\int_0^l k(x,t) \, \cdot \, dt, \quad k(x,t)=\{k_{ij}(x,t)\}_{i,j=1}^p=k_{ij}(d_i x-d_j t),
\end{equation}
where $I$ is the identity operator,
the $p \times p $ matrix function $k(x)$ on the right hand side of the
second relation in (\ref{i1}) is given by the equalities
\begin{equation}\label{i2}
k(x)=\T_2^*e^{i x \b^*}\T_1 \quad (x>0), \quad k(-x)=k(x)^*,
\end{equation}
$\T_m$ ($m=1,2$) is an $n \times p$ matrix, and $\b$ is an $n \times n$ matrix
for some integer $n>0$.
Without loss of generality we assume further that
\begin{equation}\label{i3}
d_1 \geq d_2 \geq \ldots \geq d_p>0.
\end{equation}
\begin{Rk}\label{defS}
We suppose that  equalities (\ref{i2}) hold on $(0, d_1l)$, and so, according to (\ref{i3}),
each entry $k_{ij}(x)$ is determined by  (\ref{i2}) on the interval, which contains
$(-d_j l, d_il)$, i.e., the operator $S$ of the form (\ref{i1}) is determined by (\ref{i2}).
\end{Rk}
Introduce the operator 
\begin{equation}\label{i4}
E  \in \{L^2_p(0,l), \, L^2(D)\} \quad L^2(D):=L^2(0, d_1l)\oplus L^2(0, d_2l)\oplus \ldots \oplus L^2(0, d_pl)
\end{equation}
by the equality $(Ef)_j(z)=f_j(z/d_j)$. We shall denote also by $E$ the corresponding  operator from
$\{L^2_p(0,l), \, L^2_p(0,d_1l)\}$ with the natural embedding of $L^2(D)$ into $L^2_p(0,d_1l)$:
\begin{equation}\label{i5}
(Ef)_j(z)=f_j(z/d_j) \quad (0<z<d_j l), \quad (Ef)_j(z)=0 \quad (d_jl<z<d_1 l).
\end{equation}
By  (\ref{i1}) and  (\ref{i5}) it is easy to see that
\begin{eqnarray}
\label{i6}&&
S=E^{-1}\Big(I+\int_0^a \wt k(y, z) \, \cdot \, dz\Big) E, \quad a:=d_1l, \\ 
&&
\label{i7}
 \wt k(y,z)=\{\wt k_{ij}(y, z)\}_{i,j=1}^p, \quad \wt k_{ij}(y, z)=0
\quad  {\mathrm{if}}\, \, z>d_j l \, \,
{\mathrm{or}}\, \, y> d_i l, \\
&&
\label{i8}
\wt k_{ij}(y,z)=\frac{1}{d_j}k_{ij}(y- z) \quad  {\mathrm{if}}\, \, 0<z<d_j l \, \,
{\mathrm{and}}\, \, 0<y<d_il.
\end{eqnarray}
According to (\ref{i2}), (\ref{i7}), and (\ref{i8}),  the operator 
\begin{equation}\label{i9}
\wt S=I+\int_0^a \wt k(y, z) \, \cdot \, dz
\end{equation}
is not an operator
with a difference kernel but it is a semiseparable operator.
Recall \cite{GGK1} that the integral operator $\wt S$ of the form
(\ref{i9})
 is called semiseparable, when $\wt k$ admits
representation 
\begin{equation}\label{2v3} 
\wt k(y, z)=F_1(y)G_1(z) \quad {\mathrm{for}} \, \,
 y> z, \quad \wt k(y, z)=F_2(y)G_2( z)
 \quad {\mathrm{for}} \, \, y< z,
\end{equation}
where $F_1$ and $F_2$ are $p \times n$ matrix functions and 
$G_1$ and $G_2$ are $n \times  p$ matrix functions for some $n>0$.
It is assumed that the entries of $F_1$, $F_2$, $G_1$, and $G_2$
are square integrable.
When the operator
$\wt S$ is invertible and its kernel $\wt k$ is given by (\ref{2v3}),
the kernel of the operator $\wt T=\wt S^{-1}$ is expressed in terms
of the $2 n \times 2 n$ solution $U$ of the differential equation
\begin{equation}  \label{2v4}
\Big(\frac{d}{d y}U\Big)(y)=\wt J \wt H(y)U(y), \quad y \geq 0, \quad U(0)=I_{2n}, 
\end{equation}
where
\begin{eqnarray}&&  \label{2v5}
\wt J \wt H(y):=B(y)C(y),  \quad 
\wt J=\big( \wt J^*\big)^{-1}=
\left[ \begin{array}{lr}  0 &
-I_p \\ I_p & 0
\end{array} \right].
\\  && \label{2v5'} B(y)=\left[\begin{array}{c}
-G_1(y) \\ G_2(y)
\end{array}
\right] , \quad
C(y)=\left[\begin{array}{lr}
F_1(y) & F_2(y)
\end{array}
\right].
\end{eqnarray}
Namely, we have (see, for instance, \cite{GGK1})
\begin{equation}  \label{2v6}
\wt T=\wt S^{-1}=I+\int_0^a \wt T(y,z) \, \cdot \, dz, 
\end{equation}
\begin{equation}  \label{2v7}
 \wt T(y,z)=\left\{\begin{array}{l}
C(y)U(y)\big(I_{2 n}-P^{\times}\big)U(z)^{-1}B(z), \quad y>z, \\
-C(y)U(y)P^{\times}U(z)^{-1}B(z), \quad y<z. \end{array} \right.
\end{equation}
Here $P^{\times}$ is given in terms of the $n \times n$ blocks
$U_{21}(a)$ and $U_{22}(a)$ of $U(a)$:
\begin{equation}  \label{2v8}
P^{\times}=\left[\begin{array}{lr}
0 & 0 \\ U_{22}(a)^{-1}U_{21}(a) & I_{n}
\end{array}
\right],
\end{equation}
and the invertibility of  $U_{22}(a)$ is a necessary and sufficient condition
for the invertibility of $\wt S$.

When the semiseparable operator $\wt S$ is not invertible, its kernel subspace
is given by the equality
(\cite{GGK1}, p. 157):
\begin{equation}  \label{do1}
\ker \wt S =\{h(y): \, h(y)=C(y)U(y)\left[ \begin{array}{c} 0 \\
g
\end{array} \right], \,  U_{22}(a)g=0\}.
\end{equation}

Rewrite  $D$ in the form
\begin{eqnarray} &&\nonumber
D={\mathrm{diag}}
\{\wt d_1I_{p_1},  \ldots, \wt d_k I_{p_k}  \}, \quad p_1+ \ldots +p_k=p, \\ &&
 \label{i10}
 \wt d_{j_1} > \wt  d_{j_2}> 0 \quad ( j_1<j_2 \leq k),
\end{eqnarray}
and put
\begin{equation} \label{i11}
\wt d_{k+1}=0, \quad P_{k+1}=I_p, \quad
P_j= {\mathrm{diag}}
\{I_{p_1},  \ldots,  I_{p_{j-1}}, \, 0, \ldots, 0  \} \, \, (2 \leq j \leq k).
\end{equation}
Then, in view of  (\ref{i2}), (\ref{i7}), (\ref{i8}), and (\ref{2v5'}) we have
\begin{equation}  \label{i12}
B(y)=e^{-y {\cal A}}
\left[\begin{array}{c}
- \T_1 \\ \T_2
\end{array}
\right] D^{-1}P_j, \quad
C(y)=P_j \left[\begin{array}{lr}
\T_2^* & \T_1^*
\end{array}
\right]e^{y {\cal A}},
\end{equation}
for
\begin{equation}  \label{i13}
\wt d_{j}l<y<\wt d_{j-1} l \quad (2 \leq j \leq k+1),
\quad
{\cal A}:= i\left[\begin{array}{lr}
\b^* &0 \\ 0 & \b
\end{array}
\right].
\end{equation}
\begin{Rk}\label{Ker}
By   (\ref{i5}),  (\ref{do1}) and  (\ref{i12}), it is immediate
that 
\begin{equation}  \label{do2}
\ker \, \wt S \in {\mathrm{Im}}E.
\end{equation}
where {\rm{Im}} means image.
The integral parts of $S$ and $\wt S$ are compact operators.
Hence, if $\wt S$ is not invertible, then $\ker \wt S\not=0$,
and  according to (\ref{do2}) the subspace $E^{-1}\ker \wt S$
is well defined.  In view of (\ref{i6}) and (\ref{i9}), we have $S E^{-1}\ker \wt S=0$,
i.e., $S$ is not invertible too. It follows from (\ref{i6}), (\ref{2v6}), and (\ref{2v7})  that if $\wt S$ is invertible,
then $S$ is invertible. In other words,
$S$ and $\wt S$ are simultaneousöy invertible.
\end{Rk}

Next, introduce notations
\begin{equation}  \label{i14}
 {\cal A}_j^{\times}= {\cal A}+Y_j,
 \quad
 Y_j=
\left[\begin{array}{c}
- \T_1 \\ \T_2
\end{array}
\right] D^{-1}P_j\left[\begin{array}{lr}
\T_2^* & \T_1^*
\end{array}
\right].
\end{equation}
For $2 \leq j \leq k+1$, put
\begin{equation}  \label{i15}
U(y)=e^{-y {\cal A}}e^{(y-\wt d_jl) {\cal A}_j^{\times}}e^{\wt d_jl {\cal A}}U(\wt d_jl)
\quad (\wt d_{j}l\leq y \leq \wt d_{j-1} l), \quad U(0)=I_{2n},
\end{equation}
Now, we are prepared to formulate the inversion theorem.
\begin{Tm}\label{S1}
Let $S$ be an operator with the $D$-difference kernel, which
has the form (\ref{i1}), where $k$ is given by (\ref{i2})
and $D$ satisfies (\ref{i10}). Let also  $\det U_{22}(a)\not=0$
for $U$ given by (\ref{i15}). Then $S$
is invertible and its inverse is given by the formula
$S^{-1}=E^{-1} \wt TE$, where $E$ is defined by (\ref{i5})
and $\wt T$ is given by (\ref{2v6})-(\ref{2v8}).
The  matrix functions $B$ and $C$ in (\ref{2v7})
are given by (\ref{i12}) and the $\wt J$-unitary matrix function
$U$ in (\ref{2v7}) has the form (\ref{i15}).
\end{Tm}
\begin{proof}.
To prove the theorem we need to show that $U$ of the form (\ref{i15})
satisfies  (\ref{2v4}).
Then by the properties of the semiseparable operators  we shall
obtain that $\wt S$ given by  (\ref{i9}) is invertible and that $\wt T=\wt S^{-1}$ is given by  (\ref{2v6})-(\ref{2v8}),
(\ref{i15}).  The formula $S^{-1}=E^{-1}\wt TE$ will be immediate from (\ref{i6}).

By formulas (\ref{i12}) and  (\ref{i14})
it is easy to see that $U$ of the form (\ref{i15})
satisfies equation
\begin{eqnarray} && \nonumber  
\Big(\frac{d}{d y}U\Big)(y)=
e^{-y {\cal A}}\Big( {\cal A}_j^{\times}
-{\cal A}
 \Big)
 e^{y {\cal A}}e^{-y {\cal A}}
 e^{(y-\wt d_jl) {\cal A}_j^{\times}}e^{\wt d_jl {\cal A}}U(\wt d_jl)
 \\  &&   \label{i16}
 =e^{-y {\cal A}}Y_je^{y {\cal A}}U(y)=B(y)C(y)U(y)
\end{eqnarray}
for $0 \leq y \leq a$.
Hence, by   (\ref{2v5}) $U$ satisfies (\ref{2v4}). 

Finally, let us prove that $U$ is $\wt J$-unitary, i.e., $U(y)^*\wt J U(y)=\wt J$.
Indeed, according (\ref{i13})  we have
\begin{equation}  \label{i18}
{\cal A}^*=-\wt J{\cal A}\wt J^*.
\end{equation}
As we noted in (\ref{i16}),
the equality $B(y)C(y)=e^{-y {\cal A}}Y_je^{y {\cal A}}$ is true.
Thus, taking into account (\ref{2v5}), (\ref{i14}), and (\ref{i18})
we obtain
\begin{eqnarray} && \nonumber  
\wt H(y)=\wt J^* e^{-y {\cal A}}Y_je^{y {\cal A}}=e^{y {\cal A}^*}\wt J^*Y_je^{y {\cal A}}\\  &&   \label{i19}
 =e^{y {\cal A}^*}\left[\begin{array}{c}
\T_2 \\ \T_1
\end{array}
\right] D^{-1}P_j\left[\begin{array}{lr}
\T_2^* & \T_1^*
\end{array}
\right]e^{y {\cal A}} \geq 0.
\end{eqnarray}
It follows from  (\ref{i19}) that $\wt H^*=\wt H$. Therefore, formulas (\ref{2v4}) and  (\ref{2v5})
  imply $\frac{d}{dy}\Big(U(y)^*\wt JU(y)\Big)=0$. Moreover, from $\frac{d}{dy}\Big(U(y)^*\wt JU(y)\Big)=0$
and $U(0)=I_{2n}$ we get
$U(y)^*\wt J U(y)=\wt J$.
\end{proof}
\begin{Rk}\label{T} If $S$ is invertible, then from Theorem \ref{S1} we derive
\begin{equation}  \label{i20}
T=S^{-1}=I+\int_0^l \{T_{ij}(x,t)\}_{i,j=1}^p \, \cdot \, dt, 
\end{equation}
where for $d_ix>d_jt$ and $e_i=\left[\begin{array}{lr}
\overbrace{0 \quad \ldots \quad 0}^{i-1} \quad 1 \quad 0 \quad \ldots \quad 0 
\end{array}
\right]$ we have
\[
 T_{ij}(x,t)=
e_i
\left[\begin{array}{lr}
\T_2^* & \T_1^*
\end{array}
\right]e^{d_ix {\cal A}}
U(d_ix)\big(I_{2 n}-P^{\times}\big)U(d_jt)^{-1}
e^{-d_j t {\cal A}}
\left[\begin{array}{c}
- \T_1 \\ \T_2
\end{array}
\right]e_j^*,
\]
and for  $d_i x<d_j t$ we have
\[
 T_{ij}(x,t)=
-e_i
\left[\begin{array}{lr}
\T_2^* & \T_1^*
\end{array}
\right]e^{d_ix {\cal A}}
U(d_ix)P^{\times}U(d_jt)^{-1}
e^{-d_j t {\cal A}}
\left[\begin{array}{c}
- \T_1 \\ \T_2
\end{array}
\right]e_j^*.
\]
\end{Rk}
\section{Operator identities for   operators
with $D$-difference kernels} \label{OpId}
\setcounter{equation}{0}
According to \cite{SaL1} (Ch. 6) a bounded in $L^2_p(0,l)$ operator $S$
with $D$-difference kernel, that is, an operator of the form (\ref{0.1}), (\ref{0.2}) satisfies
the operator identity
\begin{equation}\label{1.1}
AS-SA^*=i \Pi J\Pi^*,  
\end{equation}
where $A_l=A \in  \{L^2_p(0,l), \, L^2_p(0,l)\}$,  $\Pi_l=\Pi=[\Phi_1 \quad \Phi_2]$,  $\Phi_k \in
\{\BC^p, \, L^2_p(0,l)\}$, the index "$l$" is often omitted in our notations, and
\begin{equation}\label{1.2}
A=i D \int_0^x \, \cdot \, dt,  \quad \Phi_1 g =D s(x,0)g, \quad  \Phi_2 g \equiv g.
\end{equation}
It is said that $A$, $S$, and $\Pi$, which satisfy  (\ref{1.1}), form an $S$-node.
Further we assume that $A$ and $\Phi_2$ have the form (\ref{1.2}).
Operator identities play an important role in the study of structured
operators \cite{SaL1, SaL20, SaL3}.

Let us show that not only the operator with the $D$-difference kernel
satisfies (\ref{1.1})  but
the inverse statement is also true, i.e.,  (\ref{1.1}) implies that $S$ is an operator with
a $D$-difference kernel (see also the corresponding statement in Example 1.2, p. 104
\cite{SaL3}). 
Quite similar to the proof of Theorem 1.3 (\cite{SaL1}, p. 11), where
the case $D=I_p$ was treated, one can prove the following theorem
\begin{Tm}\label{TmId} Suppose a bounded operator $T \in \{L^2_p(0,l), \, L^2_p(0,l)\}$
satisfies the operator identity
\begin{equation}\label{p1}
TA-A^*T=i \int_0^l Q(x,t) \, \cdot \, dt,
\end{equation}
\begin{equation}\label{p1'}
Q(x,t)=Q_1(x)Q_2(t),
\end{equation}
where $Q$, $Q_1$, and $Q_2$ are $p \times p$,  $p \times \wh p$, and $\wh p \times p$ $(\wh p >0)$
matrix-functions, respectively. Then $T$ has the form
\begin{equation}\label{p2}
Tf=\frac{d}{dx} \int_0^l
\frac{\p}{\p t}\Upsilon(x,t)f(t)dt,  
\end{equation}
where $\Upsilon(x,t)=\{\Upsilon_{ij}(x,t)\}_{i,j=1}^p$ is absolutely continuous in $t$, and
\begin{equation}\label{p3}
 \Upsilon_{ij}(x,t):=(2d_i d_j)^{-1}\int_{d_i x+d_j t}^{f_{min}}
Q\Big(\frac{u+d_i x-d_j t}{2d_i}, \frac{u-d_i x+d_j t}{2d_j}\Big)du, 
\end{equation}
\begin{equation}
\label{p4}  f_{min}:=\min\big(d_i(2l-x)+d_jt, \, d_i x+d_j(2l-t)\big).
\end{equation}
\end{Tm}

In fact, Theorem \ref{TmId} is true for a much wider class of functions $Q$ than the one
given by (\ref{p1'}). Similar to  Theorem 2.2 in \cite{SaL1}, the next theorem is immediate
from Theorem  \ref{TmId} and equality $\wh U A \wh U =A^*$ ($(\wh U f)(x)=\ov{f(l-x)}$).
\begin{Tm}\label{TmId2}
Suppose $S \in \{L^2_p(0,l), \, L^2_p(0,l)\}$ satisfies the operator
identity $AS-SA^*=i \int_0^l\big(\Phi_1(x)+\wh \Phi_1(t)\big) \, \cdot \, dt$,
where $\Phi_1(x)$ and $\wh \Phi_1(t)$ are $p \times p$ matrix functions with the entries  from  $L^2(0,l)$.
Then $S$ is an operator with a $D$-difference kernel, i.e., the operator of the form
(\ref{0.1}), (\ref{0.2}), and $s(u,0)=D^{-1}\Phi_1(u)$,  $s(0,u)=-D^{-1}\wh \Phi_1(u)$.
Moreover, when (\ref{1.1}) holds, that is, $\wh \Phi_1(t)=\Phi_1(t)^*$
we have
\begin{equation}
\label{p12} s(x,t)=-D^{-1}s(t,x)^*D, \quad S=S^*.
\end{equation}
\end{Tm}

\section{Positive operators $S$ and an inverse problem for canonical system} \label{InvPr}
\setcounter{equation}{0}
Operators with $D$-difference kernels are essential for the construction
of   solutions of an
inverse problem for an important subclass of canonical systems   \cite{SaA1, SaL2, SaL3}.
Canonical system is a system of the form
\begin{equation}\label{1.3}
\frac{d}{d x}w(x,\lambda )=i\lambda JH(x)w(x,\lambda ),
\quad
H(x) \geq 0, \quad   J=
\left[ \begin{array}{lr}  0 &
I_p \\ I_p & 0
\end{array} \right],
\end{equation}
where the Hamiltonian $H$ is  a $m \times m$ ($m=2p$) locally summable matrix function.
A Weyl  function of the canonical system on the semi-axis $x \geq 0$ is a $p \times p$ 
 matrix function $\vp(\la)$,
which is analytic in $\BC_+$ and satisfies the condition \cite{SaL3}
\begin{equation} \label{1.6}
\int_0^\infty
\left[ \begin{array}{lr}
I_p &i \varphi (\la)^* \end{array} \right]
  w(x, \la)^*H(x)w(x, \la)
\left[ \begin{array}{c}
I_p \\ -i \varphi (\la) \end{array} \right]
dx            < \infty, \quad \la \in \BC_+.
\end{equation}
The corresponding inverse problem is the problem to recover $H$ or,
equivalently, canonical system from the Weyl function.
In the case of  rational Weyl matrix  functions several inverse problems were solved explicitly
using a GBDT version of the B\"acklund-Darboux transformation \cite{FKS1, GKS1, GKS6, MST, SaA2}.
(See \cite{D,  GeT, Gu, MS, SaA2, ZM} and references therein for various versions
of the  B\"acklund-Darboux transformation and commutation methods.)
However, taking into account that the positivity of operators $S$
and the application of the inversion formulas for semiseparable operators is of 
independent interest, we shall use a general scheme \cite{SaL0, SaL3}
and its modification \cite{SaA1}  for the inverse problem treated in this section.
As a result of the application of the general scheme to rational matrix functions,
semiseparable operators appear.
Inverse problems for self-adjoint
and skew-self-adjoint Dirac-type systems were studied 
using semiseparable operators in \cite{AGLKS}
and \cite{FKS}, respectively. 

Consider rational Herglotz $p \times p$ matrix functions $\vp$.
The statement below is immediate from Theorem 5.2 \cite{GKS6}.
\begin{Pn}\label{RkW}
If  $\vp$ is a rational matrix function such that
\begin{equation} \label{fW}
\lim_{\la \to \infty} \vp(\la)=\frac{i}{2}D, \quad \Im \vp( \la) \geq 0 \quad (\la \in \BC_+),
 \end{equation}
then $\vp$ admits a representation (i.e., realization in terms of control theory)
\begin{equation} \label{p6}
\vp(\la)=\frac{i}{2}D+\T_1^*(\b - \la I_n)^{-1}\T_2,  
\end{equation}
where $\T_1$ and $\T_2$ are $n \times p$ matrix functions, $n$ is some positive
integer number, and $n \times n$ matrix $\b$ satisfies the matrix identity
\begin{equation} \label{p6'}
 \b^*-\b =i\big(\T_2 -\T_1\big) D^{-1} \big(\T_2 -\T_1\big)^*.
\end{equation}
\end{Pn}
Dirac systems and  Weyl matrix functions $\wt \vp$, which have the form  $\wt \vp =2D^{-\frac{1}{2}}\vp D^{-\frac{1}{2}}$,
were studied in \cite{GKS1}.
The next proposition follows  from the Step 1 of the proof of Theorem 4.3 \cite{GKS1}
 (see also \cite{GKS6}).
\begin{Pn}\label{NevRepr} Let  relations  (\ref{p6}) and  (\ref{p6'}) hold. Then $\Im \vp(\la)>0$
$(\la \in \BC_+)$  and $\vp$ admits Herglotz representation
\begin{equation} \label{p7}
 \varphi ( \lambda )= \nu +  \int_{- \infty }^{ \infty }
\Big( \frac{1}{z- \lambda } - \frac{z}{1+z^{2}}\Big)d \tau (z) \quad (\nu= \nu^*),
\end{equation}
where
\begin{equation}    \label{p8}
\tau ( z )= \int_{0}^{ z }\rho(t)d t + \sum_{ z_{k}< z } \nu_{k},
\end{equation}
numbers $z_1<z_2< \ldots $ are the real eigenvalues of $\b$,
\begin{equation} \label{p9}
\nu_{k}={ \mathrm {res}}_{z=z_{k}} \T_{2}^{*}(zI_{n}- \beta )^{-1}
\T_{2} \geq 0, 
\end{equation}
and $\rho$ is $p \times p$ rational matrix function:
\begin{equation} \label{p10}
\rho( t)=\frac{1}{2 \pi}\zeta( t)^* D\zeta( t) \geq 0, \hspace{1em}
\zeta( t):=
I_{p}-i D^{-1} \big(\T_2 -\T_1\big)^*
( t I_{n}-   \beta )^{-1}
 \T_{2}.
 \end{equation}
\end{Pn}

It is easy to see from  (\ref{1.2}) that
\begin{equation} \label{p11}
(I-zA)^{-1}\Phi_2=e^{izxD}, \quad \Phi_2^*(I-zA^*)^{-1}f=\int_0^le^{-izxD}f(x)dx \quad (f \in L^2_p(0,l)).
\end{equation}
By (\ref{p8})-(\ref{p11}) the right-hand side of the equality
\begin{equation} \label{p13}
S:=\int_{-\infty}^{\infty}(I-zA)^{-1}\Phi_2d \tau(z)\Phi_2^*(I-zA^*)^{-1}
\end{equation}
weakly converges, and so the equality defines an operator $S$. Moreover, it is easy to see that
the inequalities  
\begin{equation} \label{p13'}
c(f,f)_{L^2}>(Sf,f)_{L^2}>0
\end{equation}
 hold for some fixed $c>0$ and arbitrary  $f \not=0$.
(Here $(\cdot , \cdot)_{L^2}$ denotes the scalar product in $L^2_p(0,l)$.)
 Thus,
$S$ is a bounded and positive operator. We shall show that operators $S$ 
belong to a subclass of operators of the form (\ref{i1}),  (\ref{i2}).
\begin{Tm}\label{TmPos}
  Let the matrix identity  (\ref{p6'})
hold. Then the operator $S$ given by   (\ref{i1})  and  (\ref{i2}) is positive
 and boundedly invertible.
 \end{Tm}
\begin{proof}. The theorem is obtained by proving that $S$ of the form
(\ref{i1}),  (\ref{i2}) admits representation (\ref{p13}).

First, consider  $S$ given by  (\ref{p13}).
It can be calculated directly (see also Section 1.1 in \cite{SaL20}) that this operator
$S$  satisfies the operator identity  (\ref{1.1}), where $\Pi=[\Phi_1 \quad
\Phi_2]$ and
\begin{equation} \label{p14}
\Phi_1=i\left(\nu -\int_{-\infty}^{\infty}\Big(A(I-zA)^{-1}+\frac{z}{1+z^2}I\Big)\Phi_2d \tau(z)
\right).
\end{equation}
Here the operator $\Phi_1$ is an operator of multiplication by the matrix function,
which we denote by $\Phi_1(x)$. From the identity (\ref{1.1}) and Theorem \ref{TmId2} it follows
that $S$ is an operator with a $D$-difference kernel 
$s(x,t)=\{ s_{ij}(d_ix-d_j t)\}_{i,j=1}^p$ and  $s(x,0)=D^{-1}\Phi_1(x)$.
Introduce $S=S_l$ and $\Phi_1=\Phi_{1,l}$ by (\ref{p13}) and (\ref{p14}), respectively, for all $0<l<\infty$. 
Then the kernel $s(x)$
of the integral operators $S_l$ is determined on $\BR$ by the equalities
\begin{equation} \label{p15}
s_{ij}(x)=d_i^{-1}\big(\Phi_1\big)_{ij}(x/d_i) \quad (x>0), \quad s_{ij}(-x)=-\frac{d_j}{d_i}\ov{s_{ji}(x)}.
\end{equation}
For $\vp$ satisfying  (\ref{p7}),  according to Statement 3  in \cite{SaA1}, 
after the corresponding change of notations we get
\begin{equation} \label{1.7}
\varphi (\la)=\la \int_0^{\infty}s(x,0)^*e^{i\la x D}dx D^2=\la \int_0^{\infty}e^{i\la x }s(x)^*dx D.
\end{equation}
Note that in view of    formula (\ref{p14}) and Proposition \ref{NevRepr} we can present $s$ as a sum
$s(x)=s_1(x)+s_2(x)$, where the entries of $s_1$ are bounded and the entries of $s_2$
belong $L^2(0,\infty)$. Finally, we apply Fourier transform to derive from (\ref{1.7})
the equality
\begin{equation} \label{1.7'}
e^{-\eta x}s(x)^*
=\frac{1}{2 \pi}
{\mathrm{
l.i.m.}}_{a \to \infty} \int_{- a}^{a}e^{-i \xi
x} \lambda^{-1} \vp(\lambda ) D^{-1}d \xi \quad (\lambda= \xi +i \eta , \quad \eta>0),
\end{equation}
the limit l.i.m. being the limit in $L^2(0,l)$ ($0<l<\infty$). Using  (\ref{p6}) and (\ref{1.7'}), we obtain
\begin{equation} \label{p16}
e^{-\eta x}s(x)^*
=\frac{1}{2 \pi}
{\mathrm{
l.i.m.}}_{a \to \infty} \int_{\G_a}e^{-i \xi
x} \lambda^{-1} \vp(\lambda ) D^{-1}d \xi \quad (\lambda= \xi +i \eta , \quad \eta>0),
\end{equation}
where $\G_a$  is a clockwise oriented contour: 
\[
\G_a=[-a, \, a]\cup\{\xi: \, |\xi|=a, \, \Im \xi<0\}.
\]
It is easy to see that
\begin{equation} \label{p16'}
\frac{1}{2 \pi}
{\mathrm{
l.i.m.}}_{a \to \infty} \int_{\G_a}e^{-i \xi
x} \lambda^{-1} d \xi =-ie^{-\eta x}.
\end{equation}
According to  (\ref{p6'}) we have $\s(\b)\subset \ov{\BC_-}$,
where $\s$ is spectrum.
Similar to \cite{FKS} we turn to zero $\ve$ in the equality
$\la^{-1}(\b_{\ve}-\la I_n)^{-1}=\b_{\ve}^{-1}\big(\la^{-1}I_n+(\b_{\ve}-\la I_n)^{-1}\big)$,
where $\det \b_{\ve}\not=0$, $\|\b -\b_{\ve}\|<\ve$, and thus obtain
\begin{equation} \label{p17}
\frac{1}{2 \pi}
{\mathrm{
l.i.m.}}_{a \to \infty} \int_{\G_a}e^{-i \xi
x} \lambda^{-1} (\b-\la I_n)^{-1}d \xi =e^{-\eta x}\int_0^x\exp(-iu\b)du.
\end{equation}
Here we take into account that, when the spectrum of some matrix
${\cal{K}}$ is situated inside the
anti-clockwise oriented contour $\G$ we have
\[
\frac{1}{2\pi i}\int_{\G}e^{- i \la x}(\la I_n - {\cal{K}})^{-1}d\la=\exp(-ix{\cal K}).
\]
By (\ref{p6}) and (\ref{p16})-(\ref{p17}) we get
\begin{equation} \label{p18}
s(x)
=\frac{1}{2 }I_p+D^{-1}\T_2^*\int_0^x\exp(iu\b^*)du\T_1 \quad (x>0).
\end{equation}
It follows from  (\ref{p15}) that $s(x)=-D^{-1}s(-x)^*D$ ($x<0$), and so according
to (\ref{p18}) $s(x)$ is continuously differentiable for $x \not=0$.
As  the functions $s_{ij}(x)$ are continuous at $x=0$ for $i\not=j$,
and $s_{ii}(+0)-s_{ii}(-0)=1$, formulas  (\ref{0.1}) and  (\ref{0.2})  imply
 (\ref{i1}), where $k(x)=D\Big(\frac{d}{dx}s\Big)(x)$. Therefore we have
\begin{equation} \label{p19}
k(x)=
\T_2^*\exp(ix\b^*)\T_1 \quad (x>0), \quad k(x)=k(-x)^*.
\end{equation}

Now, note that equalities 
 (\ref{i2}) and  (\ref{p19}) coincide. In other words,
the operator $S$, which is considered in the theorem,
admits representation (\ref{p13}). Hence, by (\ref{p13'}) 
this operator is bounded and positive, and so in view of 
(\ref{i1}) and (\ref{i2}) it is also boundedly invertible.
\end{proof}
The matrix function $\tau$ of the form (\ref{p8})-(\ref{p10})
and the $S$-node given by   (\ref{1.2}), (\ref{p13}), and (\ref{p14})
satisfy conditions of Theorem 2.4 \cite{SaL3}, p. 57. 
Therefore $\vp(\la)$ given by  (\ref{p7}) can be presented as a linear-fractional transformation
\begin{equation}\label{p20}
\varphi (\la ) =i\big( {\cal W}_{11}(\la )R_1(\la )+ {\cal W}_{12}(\la )R_2(\la )\big) 
\big( {\cal W}_{21}(\la )R_1(\la )+ {\cal W}_{22}(\la )R_2(\la )\big)^{-1},
\end{equation}
where ${\cal W}_{ij}(\la )$ are $p \times p$ blocks of the matrix function ${\cal W}$,
\begin{equation}\label{p21}
{\cal W}(\la ):=W(l, \ov \la)^*, \quad W(l,  \la)=I_{2p}+i\la J\Pi^*S^{-1}(I -\la A)^{-1}\Pi,
\end{equation}
and $R_1(\la)=R_1(l,\la)$, $R_2(\la)=R_2(l,\la)$ is a pair of $p \times p$ matrix functions,
which are meromorphic in $\BC_+$
 and have property-$J$, that is,
\begin{equation}\label{p22}
R_1(\la )^*R_1(\la )+R_2(\la )^*R_2(\la )>0,\quad \left[\begin{array}{lr}
  R_1(\la )^* & R_2(\la )^* \end{array}\right]\, J\,
\left[\begin{array}{c} R_1(\la ) \\  R_2(\la ) \end{array}\right]\geq 0.
\end{equation}
It is easy to see from (\ref{p21})  that $\lim_{l \to +0}W(l,  \la)=I_{2p}$ and thus we put $W(0,\la)=I_{2p}$.
Now, by Theorem 2.1 from  \cite{SaL3}, p.54 the matrix function $W$  satisfies for $x \geq 0$
the equation
\begin{equation} \label{p23}
W(x, \la)=I_{2p}+i\la J \int_0^x\big(dB_1(r))W(r, \la), \quad B_1(r):=\Pi_r^*S_r^{-1}\Pi_r,
\end{equation}
where $S_r \in \{L^2_p(0,r), \, L^2_p(0,r)\}$,  $\Pi_r \in \{\BC^{2p}, \, L^2_p(0,r)\}$.
As the operators $S_{r}$ ($0<r \leq l <\infty$) are invertible, the operators
$S_l$ admit triangular factorisation (see \cite{GoKrb}, p. 184). It follows that $B_1$
is differentiable, and we rewrite (\ref{p23})  as the canonical system
\begin{eqnarray} \label{p24}&&
\frac{d}{d x}W(x, \la)=i\la J H(x)W(x, \la),  
\\ &&\label{p24'}
H(x):=\frac{d}{dx}\Big(\Pi_x^*S_x^{-1}\Pi_x\Big).
\end{eqnarray}
Moreover, in view of Remark \ref{T} the kernel $T_r(x,t)$ of the integral operator $S_r^{-1}$ 
is continuous  with respect to $x,t,r$ excluding the lines $d_i x = d_j t$.  Therefore, for $d_i x \not= d_j r$  ($1\leq i,j \leq p$)
similar to the continuous kernels
(\cite{GoKrb}, p.186)  we have
\begin{equation} \label{p31}
k(x,r)+T_r(x,r)+\int_0^rk(x,u)T_r(u,r)du=0, \quad x \leq r \leq l.
\end{equation}
Introduce an upper triangular operator 
\begin{equation} \label{p32}
V_+=I+\int_x^lT_r(x,r) \, \cdot \, dr \in \{L^2_p(0, \, l)\}.
\end{equation}
According to (\ref{p31}) and (\ref{p32}) the operator $S_lV_+$ is a lower  triangular operator.
Hence, the operator $V_+^*S_lV_+$ is a lower triangular operator. On the other hand
$V_+^*S_lV_+$ is selfadjoint, and so the integral part of $V_+^*S_lV_+$ equals zero, i.e.,
$V_+^*S_lV_+=I$ or equivalently
\begin{equation} \label{p33}
S_l^{-1}=V_+V_+^*, \quad V_{+,l}^*=V_+^*=I+\int_0^xT_x(x,r) \, \cdot \, dr.
\end{equation}
In the second equality above we used formula (\ref{p32}) and relation
$T_x(r,x)^*=T_x(x,r)$ ($x \geq r$).
\begin{Tm}\label{aInvPr} Let $\vp$ be a rational  function, which satisfies
(\ref{fW}). Then $\vp$ is a Weyl function of the canonical system
(\ref{p24}), where the Hamiltonian $H$
has the form
\begin{equation} \label{p34}
H(x)=\g(x)^*\g(x), \quad \g(x)=\Big(V_+^*[\Phi_1 \quad \Phi_2]\Big)(x) \quad (x\leq l<\infty),
\end{equation}
and the operator $V_+^*$ is given by (\ref{p33}) and is
applied columnwise to the matrix functions $\Phi_1(x)=\{d_is_{ij}(d_i x)\}_{i,j=1}^p$
and $\Phi_2 \equiv I_p$. The matrix function $s(x)$ is given by (\ref{p18}) and the
matrix function $T_x(x,r)$ in (\ref{p33}) is given in Remark \ref{T}.
\end{Tm}
\begin{proof}. 
It follows from  (\ref{p24})  that
\begin{eqnarray} && \label{p25}
\frac{d}{dx}\Big(W(x, \ov \la)^*JW(x,  \la)\Big)=0, \\ 
&& \label{p26}
 \frac{d}{dx}\Big(W(x,  \la)^*JW(x,  \la)\Big)
=i( \la - \ov  \la)W(x,  \la)^*H(x)W(x,  \la).
\end{eqnarray}
In view of  (\ref{p26}) we obtain
 \begin{equation} \label{p27}
\int_0^l
  W(x, \la)^*H(x)W(x, \la)
dx     =  i(\ov \la -  \la)^{-1}\Big(W(l,  \la)^*JW(l,  \la)-J\Big).
\end{equation}
Note also that according to (\ref{p25}) the equality $W(l, \ov \la)^*JW(l,  \la)=J$ holds,
or equivalently 
\begin{equation} \label{p28}
W(l, \ov \la)^*=JW(l,  \la)^{-1}J.
\end{equation}

By Proposition  \ref{RkW} $\vp$ admits representation (\ref{p6})
and identity (\ref{p6'}) is valid. So, by Proposition \ref{NevRepr} $\vp$ admits
Herglotz representation, where the matrix function $\tau(t)$
has the form (\ref{p8})-(\ref{p10}). Hence,
as it  was
shown above, the representation (\ref{p20}) of $\vp$, where ${\cal W}$
is expressed via the matrizant $W(l, \la)$ and the pair $R_1$, $R_2$ satisfies (\ref{p22}), is also true.
Using (\ref{p28}), we rewrite (\ref{p20}) in the form 
\begin{equation} \label{p29}
\left[\begin{array}{c}
I_p \\ -i\vp(\la)
\end{array}
\right]=W(l,  \la)^{-1}J\left[\begin{array}{c}
R_1(\la) \\  R_2(\la)
\end{array}
\right]\big( {\cal W}_{21}(\la )R_1(\la )+ {\cal W}_{22}(\la )R_2(\la )\big)^{-1}.
\end{equation}
Taking into account (\ref{p22}), (\ref{p27}), and (\ref{p29}) we derive
\begin{eqnarray} && \nonumber
\int_0^l
\left[ \begin{array}{lr}
I_p &i \varphi (\la)^* \end{array} \right]
  W(x, \la)^*H(x)W(x, \la)
\left[ \begin{array}{c}
I_p \\ -i \varphi (\la) \end{array} \right]
dx  \leq i( \la - \ov  \la)^{-1}\\ \label{p30} &&
\times \left[ \begin{array}{lr}
I_p &i \varphi (\la)^* \end{array} \right]
J
\left[ \begin{array}{c}
I_p \\ -i \varphi (\la) \end{array} \right],
\quad \la \in \BC_+.
\end{eqnarray}
As the right-hand side in  the inequality (\ref{p30}) does not depend on $l$
we can substitute $\infty$ instead of the limit $l$ of integration in the left-hand side.
Hence $\vp$ is a Weyl function of the constructed system.

According to the second relation in (\ref{p33}) we obtain $(V_{+,l}^*f)(x)=(V_{+,x}^*\wt f)(x)$
for $x \leq l$, where $\wt f$ is the restriction of $f$ on the interval $[0,x]$.
Therefore, relations (\ref{p24'}) and (\ref{p33}) imply (\ref{p34}).
\end{proof}
\begin{Cy} \label{appl} Let the conditions of Theorem \ref{aInvPr} hold and let $\det \b \not=0$. Then we have
\begin{equation} \label{p35}
 \g(x)=\Big(V_+^*[\frac{1}{2}D+i\T_2^*(\b^*)^{-1}\T_1  \qquad I_p]\Big)(x)-i[\g_0(x) \quad 0],
\end{equation}
where the $s$-th row of $\g_0$ $(p \geq s \geq 1)$ is given by the equality
\begin{eqnarray} \nonumber &&
e_s\g_0(x)=e_s \Big(\T_2^* e^{id_sx {\b^*}}+
\left[\begin{array}{lr}
\T_2^* & \T_1^*
\end{array}
\right]e^{d_sx {\cal A}}
U(d_s x)
\\ \label{p36} && \times
\big(P^{\times}U(d_1x)^{-1}- U(d_sx)^{-1}+ I_{2 n}-P^{\times}\big)\left[\begin{array}{c}
I_p \\ 0
\end{array}
\right]\Big)(\b^*)^{-1}\T_1,
\end{eqnarray}
$U$  in (\ref{p36}) is defined by  (\ref{i15})
 after substitution  $l=x$, and $P^{\times}$ is defined by
(\ref{2v8}) after substitution  $a=d_1l=d_1x$.
\end{Cy}
\begin{proof}.
By (\ref{0.2}), (\ref{1.2}), and (\ref{p18})  the equality
\begin{equation} \label{p37}
e_s[\Phi_1(x) \quad \Phi_2]=e_s[\frac{1}{2}D+i\T_2^*(\b^*)^{-1}\T_1 -i \T_2^* e^{id_sx {\b^*}}
(\b^*)^{-1}\T_1 \qquad I_p]
\end{equation}
is true. Using  (\ref{p37})  and the second equlity in (\ref{p34}) we obtain (\ref{p35}), where
\begin{equation} \label{p38}
\g_0(x)=\Big(V_+^*\{e_s\T_2^* e^{id_sx {\b^*}}
(\b^*)^{-1}\T_1\}_{s=1}^p\Big)(x).
\end{equation}
From (\ref{i13}) it follows that
\begin{equation} \label{p39}
\T_2^* e^{id_sx {\b^*}}=[\T_2^* \quad \T_1^*]e^{d_sx {\cal A}}\left[\begin{array}{c}
I_p \\ 0
\end{array}
\right].
\end{equation}
According to the representation of $V_+^*$ in (\ref{p33}), Remark \ref{T}, formula (\ref{p39})
and second relation in (\ref{i14}) we get
\begin{equation} \label{p40}
V_+^*\{e_s\T_2^* e^{id_sx {\b^*}}
\}_{s=1}^p=\{e_s\T_2^* e^{id_sx {\b^*}}
\}_{s=1}^p
+\{{\cal F}_s(x)
{\cal G}_s(x)
\}_{s=1}^p,
\end{equation}
where
\begin{equation} \label{p41}
{\cal F}_s(x)=e_s\left[\begin{array}{lr}
\T_2^* & \T_1^*
\end{array}
\right]e^{d_sx {\cal A}}
U(d_s x),
\end{equation}
\begin{eqnarray}&& \nonumber
{\cal G}_s(x)=
\Big((I_{2 n}-P^{\times})\int_0^{d_s x}U(z)^{-1}e^{-z{\cal A}}Y(z)e^{z {\cal A}}dz
\\ \label{p42} &&
-P^{\times}\int_{d_s x}^{d_1 x}U(z)^{-1}e^{-z{\cal A}}Y(z)e^{z {\cal A}}dz\Big)\left[\begin{array}{c}
I_p \\ 0
\end{array}
\right],
\end{eqnarray}
\[
Y(z)=\sum_{j: \, d_j >\wt d_m}\frac{1}{d_j}\left[\begin{array}{c}
- \T_1 \\ \T_2
\end{array}
\right] e_j^*e_j\left[\begin{array}{lr}
\T_2^* & \T_1^*
\end{array}
\right]=Y_m \quad {\mathrm{for}} \quad  \wt d_m x\leq z \leq \wt d_{m-1}x.
\]
Taking into account  (\ref{i16}) rewrite (\ref{p42}) in the form
\begin{equation} \label{p43}
{\cal G}_s(x)=
\Big((I_{2 n}-P^{\times})\big(I_{2 n}-U(d_s x)^{-1}\big)
+P^{\times}\big(U(d_1x)^{-1}-U(d_s x)^{-1}\big)\Big)\left[\begin{array}{c}
I_p \\ 0
\end{array}
\right].
\end{equation}
Finally, formulas  (\ref{p38}), (\ref{p40}), (\ref{p41}), and  (\ref{p43}) imply  (\ref{p36}).
\end{proof}

In view of Corollary \ref{appl}, to recover $\g$ and Hamiltonian $H$
 we need only to calculate the action of $V_+^*$  on constant
vectors.

The  matrix function
$\g(x)$, which is recovered in Theorem \ref{aInvPr}, satisfies the equality
\begin{equation} \label{p44}
\g(x)J\g(x)^*\equiv D.
\end{equation}
Indeed, by  (\ref{1.1}), the first equality in  (\ref{p33}),  and the second equality in (\ref{p34})
we have
\begin{equation} \label{p45}
V_+^*A(V_+^*)^{-1}-V_+^{-1}A^*V_+=i\g(x)J\int_0^l \g(t)^*\, \cdot \, dt.
\end{equation}
As $V_+^*A(V_+^*)^{-1}$ is a lower triangular operator and
$V_+^{-1}A^*V_+$ is an upper triangular operator, we derive
\begin{equation} \label{p46}
V_+^*A(V_+^*)^{-1}=i\g(x)J\int_0^x \g(t)^*\, \cdot \, dt.
\end{equation}
Rewrite (\ref{p46}) in terms of the kernels of the corresponding integral operators
and put $t=x$ 
to get (\ref{p44}).

As it is stated in the proposition below, equality (\ref{p44}) means that we recover canonical systems from the subclass 
of systems with
linear similar matrix functions $JH(x)$, though (differently from \cite{SaL3}, p. 104)
the kernel of $S^{-1}$ is not necessarily continuous.
\begin{Pn} \label{sim}Let the conditions of Theorem \ref{aInvPr} hold.
Then $JH(x)$ is similar to the matrix $JH_0$, where
\begin{equation} \label{p47}
H_0:=\left[\begin{array}{lr}
D & 0 \\ 0 & 0
\end{array}
\right]
\end{equation}
\end{Pn}
\begin{proof}.  Fix $x \geq 0$ and denote by $X$ a $p \times 2p$  matrix
such that  it has rank $p$ and satisfies the equality $XJ\g(x)^*=0$.  As the maximal $J$-nonnegative subspaces
are $p$-dimensional,
it  easily follows from $\g(x)J\g(x)^*>0$ and
$XJ\g(x)^*=0$ that $XJX^*<0$. Then, we have
\begin{equation} \label{p48}
\wt X J \wt X^*=-I_p, \quad \wt X J \g(x)^*=0 \quad {\mathrm{for}} \quad
\wt X:=(-XJX^*)^{-\frac{1}{2}}X.
\end{equation}
Now, put
\begin{equation} \label{p49}
L:=\left[\begin{array}{c}
D^{-\frac{1}{2}}\g(x) \\ \wt X
\end{array}
\right].
\end{equation}
By  (\ref{p44}), (\ref{p48}), and (\ref{p49})   the equality
\begin{equation} \label{p50}
L^{-1}=[J\g(x)^*D^{-\frac{1}{2}} \quad -J\wt X^*]
\end{equation}
is true. According to  (\ref{p47}), (\ref{p49}), and (\ref{p50}) we get
$L^{-1}H_0L=J\g(x)^*\g(x)$. In view of  (\ref{p34}) the last equality
yields $L^{-1}H_0L=JH(x)$.
\end{proof}
{\bf Acknowledgement.}
The work of A.L. Sakhnovich was supported by the Austrian Science Fund (FWF) under
Grant  no. Y330, and his visit to Mexico was supported by the PIFI grant 
P/CA-9 2007-14-17.
 A.L. Sakhnovich is grateful to  the Autonomous University of Hidalgo
for its hospitality.

\begin{flushright} \it
A.L. Sakhnovich, \\  Fakult\"at f\"ur Mathematik,
Universit\"at Wien,
\\
Nordbergstrasse 15, A-1090 Wien, Austria \\ 
e-mail: {\tt al$_-$sakhnov@yahoo.com }
\end{flushright}

\begin{flushright} \it
A.A. Karelin, \\
Universidad Autonoma del Estado de Hidalgo, 
Instituto de Ciencias Basicas e Ingenieria
Centro de investigaci—n Avanzada en Ingenier'a Industrial
Carretera Pachuca-Tulancingo, Km. 4,5 Ciudad Universitaria,
C.P. 42180, Pachuca, Hidalgo, Mexico
\\
karelin@uaeh.edu.mx
\end{flushright}

\begin{flushright} \it
J. Seck-Tuoh-Mora,  \\
Universidad Autonoma del Estado de Hidalgo, 
Instituto de Ciencias Basicas e Ingenieria
Centro de investigaci—n Avanzada en Ingenier'a Industrial
Carretera Pachuca-Tulancingo, Km. 4,5 Ciudad Universitaria,
C.P. 42180, Pachuca, Hidalgo, Mexico
\\
jseck@uaeh.edu.mx
\end{flushright}

\begin{flushright} \it
G. Perez-Lechuga,   \\
Universidad Autonoma del Estado de Hidalgo, 
Instituto de Ciencias Basicas e Ingenieria
Centro de investigaci—n Avanzada en Ingenier'a Industrial
Carretera Pachuca-Tulancingo, Km. 4,5 Ciudad Universitaria,
C.P. 42180, Pachuca, Hidalgo, Mexico
\\
glechuga2004@hotmail.com
\end{flushright}
\begin{flushright} \it
 M. Gonzalez-Hernandez, \\
Universidad Autonoma del Estado de Hidalgo, 
Instituto de Ciencias Basicas e Ingenieria
Centro de investigaci—n Avanzada en Ingenier'a Industrial
Carretera Pachuca-Tulancingo, Km. 4,5 Ciudad Universitaria,
C.P. 42180, Pachuca, Hidalgo, Mexico
\\
mghdez@uaeh.edu.mx
\end{flushright}


\begin{thebibliography}{AGKS}
\bibitem{AGLKS}
{D. Alpay, I. Gohberg, L. Lerer,  M.A. Kaashoek,  and A.L. Sakhnovich},
{\it Krein systems},  in: OT: Adv. Appl. {\bf 191}, 2009, 19--36. 

\bibitem{BGK}
H. Bart, I. Gohberg, and M.A. Kaashoek, {\it Convolution equations and linear systems},
IEOT {\bf 5} (1982), 283--340.

\bibitem{D}
P.A. Deift, {\it Applications of a commutation formula}, { Duke
Math. J.} {\bf 45} (1978), 267--310.

\bibitem{CF}
M.J. Corless and A.E. Frazho,  {\it Linear Systems and Control - An
Operator Perspective}, Marcel Dekker, New York,  2003.

\bibitem{FKS1}
B. Fritzsche, B. Kirstein, and A.L. Sakhnovich, 
{\it Completion problems and scattering problems for Dirac type
differential equations with singularities},  J. Math. Anal. Appl.
{\bf 317} (2006), 510--525.

\bibitem{FKS}
B. Fritzsche, B. Kirstein, and 
A.L. Sakhnovich, {\it Semiseparable integral operators  and explicit solution
of an inverse problem for the  skew-self-adjoint
Dirac type system},  arXiv:0904.2357

\bibitem{GeT}
F. Gesztesy and G. Teschl,  {\it On the double commutation
method}, Proc. Am. Math. Soc. {\bf 124}:6  (1996), 1831-1840.

\bibitem{GGK1} 
I. Gohberg, S. Goldberg, and M.A. Kaashoek, {\it Classes
of Linear Operators}, Volume I, Birkh\"auser Verlag, Basel, 1990.

\bibitem{GK84}
I. Gohberg and M.A. Kaashoek, 
{\it Time varying linear systems with boundary conditions and integral operators. 
I. The transfer operator and its properties}, 
IEOT  {\bf 7}  (1984),   325--391.

\bibitem{GKS1}
I. Gohberg, M.A. Kaashoek, and A.L. Sakhnovich, {\it Canonical systems
with rational spectral densities: explicit formulas and
applications}, { Math. Nachr.} {\bf 194} (1998), 93--125.

\bibitem{GKS6}
I. Gohberg, M.A. Kaashoek, and A.L. Sakhnovich, {\it Scattering
problems for a canonical system with a pseudo-exponential
potential}, Asymptotic Analysis, {\bf 29}:1 (2002), 1--38.



\bibitem{GoKaSc}
I. Gohberg, M.A.  Kaashoek and F. van Schagen, 
{\it On inversion of convolution integral operators on a finite interval,}  
OT: Adv. Appl. {\bf 147} (2004), Birkh\"auser, Basel, 277--285.

\bibitem{GoKr}
I.C. Gohberg and M.G. Krein,  {\it Systems of integral equations on a half line with kernels 
depending on the difference of arguments},  Amer. Math. Soc. Transl. (2)  {\bf 14}  (1960),  217--287.

\bibitem{GoKrb}
I.Gohberg  and  M.G.Krein,  {\it Theory  and  applications  of
Volterra operators  in  Hilbert  space},  Nauka, Moscow,  1967.
Translated  in:  Transl.  of  math.  monographs  {\bf 24},
Providence,  Rhode Island,  1970.

\bibitem{Gu}
C.H. Gu, H. Hu, and Z. Zhou, {\it  Darboux transformations in
integrable systems}, Springer Verlag, 2005.



\bibitem{Ka}
M. Kac, 
{\it On some connections between probability theory and differential and integral equations,} 
Proc. Berkeley Sympos. Math. Statist. Probability, California Juli 31-- August 12, 1950 (1951), 189--215.

\bibitem{KarLT}
A. A. Karelin, Kh.  Peres Lechuga, and A.A. Tarasenko, {\it The Riemann problem and singular 
integral equations with coefficients generated by piecewise-constant functions}. 
(Russian)  Differ. Uravn.  {\bf 44}:9  (2008),  1182--1192.

\bibitem{Kr}
M.G. Krein,  {\it 
Integral equations on the half-line with a kernel depending on 
the difference of the arguments},  (Russian)
Uspehi Mat. Nauk {\bf 13}:5(83) (1958), 3--120. 

\bibitem{MS}
V.B. Matveev and M.A. Salle, {\it Darboux transformations and
solitons}, Springer Verlag, Berlin, 1991.

\bibitem{MST}
R. Mennicken, A.L. Sakhnovich, and C. Tretter, {\it Direct and
inverse spectral problem for a system of differential equations
depending rationally on the spectral parameter},  Duke Math. J.
{\bf 109}:3 (2001), 413--449.


\bibitem{SaA1}
A.L. Sakhnovich, {\it Asymptotics of spectral functions of an
$S$-node},  Soviet Math. (Iz. VUZ) {\bf 32} (1988), 92--105.

\bibitem{SaA2}
A.L. Sakhnovich, {\it Iterated B\"acklund-Darboux transform for
canonical systems}, { J. Functional Anal.} {\bf 144} (1997),
359--370.

\bibitem{SaA3} A.L. Sakhnovich, {\it Generalized
B\"acklund-Darboux transformation: spectral properties and
nonlinear equations}, JMAA {\bf 262} (2001), 274-306.

\bibitem{SaLdiff}
L.A. Sakhnovich, {\it Equations with a difference kernel on a finite interval},
{Russian Math. Surv.} {\bf 35}  (1980), 81--152.

\bibitem{SaL0}
L.A. Sakhnovich, {\it Factorisation  problems  and  operator
identities}, { Uspekhi  Mat. Nauk}  { \bf 41  }:1 (1986),  3--55;
English transl.  in {  Russian
 Math. Surveys} { \bf 41} (1986), 1-64.

\bibitem{SaL1}
L.A. Sakhnovich, {\it Integral equations with difference kernels
on finite intervals}, Operator Theory: Adv. Appl. {\bf 84},
Birkh\"auser, Basel-Boston-Berlin, 1996.

\bibitem{SaL20}
L.A. Sakhnovich, {\it Interpolation theory and its applications,} 
Mathematics and its Applications {\bf 428}, Kluwer Academic
Publishers, Dordrecht, 1997.

\bibitem{SaL2}
L.A. Sakhnovich, 
{\it On a class of canonical systems on half-axis,}
IEOT {\bf 31} (1998),  92-112.

\bibitem{SaL3}
L.A. Sakhnovich,  {\it Spectral theory of canonical differential
systems, method of operator identities}, OT: Adv.
Appl. {\bf 107}, Birkh\"auser Verlag,  Basel-Boston, 1999.

\bibitem{VVGM}
R. Vandebril, M.  Van Barel, G.  Golub, and N. Mastronardi,  
{\it A bibliography on semiseparable matrices,}
Calcolo  {\bf 42}  (2005),   249--270.

\bibitem{ZM}
V.E.Zakharov and A.V.Mikhailov, {\it On the integrability of
classical spinor models in two-dimensional space-time}, Comm.
Math. Phys. {\bf 74} (1980), 21--40.





\end{thebibliography}
\end{document}